\documentclass[12pt,a4paper]{article}
\usepackage{amsmath,amssymb,amsfonts}
\def\Bbb{\mathbb}

\title{\bf The Period Length of Euler's Number $e$}

\author{Kurt Girstmair}
\date{}

\makeatletter
\let\@@maketitle=\maketitle
\def\maketitle{\def\thispagestyle##1{\relax}\@@maketitle}
\makeatother
\newtheorem{theorem}{Theorem}

\newtheorem{lemma}{Lemma}

\def\BE{\begin{equation}}
\def\EE{\end{equation}}
\def\BD{\begin{displaymath}}
\def\ED{\end{displaymath}}
\def\BA{\begin{array}}
\def\EA{\end{array}}
\def\BEA{\begin{eqnarray*}}
\def\EEA{\end{eqnarray*}}
\def\BI{\bibitem}

\def\Z{\Bbb Z}
\def\Q{\Bbb Q}
\def\R{\Bbb R}

\def\phi{\varphi}
\def\EPS{\varepsilon}

\def\CMOD#1#2#3{#1 \equiv #2  \mbox{ mod } #3}

\def\MB{\mbox}
\def\LD{\ldots}
\def\OV{\overline}

\def\BQ{``}

\def\EQP{''}

\def\MN{\medskip\noindent}
\def\STOP{\hfill$\blacksquare$}

\def\JS#1#2{ \left( \frac{#1}{#2} \right) }

\def\END{\\ \hline\rule{0mm}{5mm}}

\begin{document}
\maketitle

\begin{abstract}

\noindent
Let $s_k/t_k$, $k\ge 0$, be the convergents of the continued fraction expansion of a number $x\in\R\smallsetminus\Q$.
We investigate the sequence of Jacobi symbols $\JS{s_k}{t_k}$, $k\ge 0$. We show that this sequence is
purely periodic with shortest possible period length 24 for $x=e=2.718281\LD$ and shortest possible
period length 40 for $x=e^2$. Further, we make
the first steps towards a general theory of such sequences of Jacobi symbols. For instance, we show that there
are uncountably many numbers $x$ such that this sequence has the period 1 (of length 1), and that every natural number
$L$ actually occurs as the shortest possible period length of some $x$.

\end{abstract}

\section*{Introduction}

Let $(a_0,a_1, a_2,\LD)$ be the {\em regular continued fraction expansion} of $x\in\R\smallsetminus\Q$.
The sequence
\BD
\label{0.1}
   s_k/t_k=[a_0,\LD,a_k], k\ge 0,
\ED
of convergents of $x$ is defined in the well-known way by
\BE
\label{0.2}
  \begin{array}{lll}
   s_{-1}=1,& s_0=a_0,& s_k=a_ks_{k-1}+s_{k-2},\\
   t_{-1}=0,& t_0=1,& t_k=a_kt_{k-1}+t_{k-2},\enspace k\ge 1.
  \end{array}
\EE
Accordingly, we may write
\BE
\label{0.3}
   x=\lim_{k\to\infty}[a_0,\LD,a_k]=[a_0,a_1,a_2,\LD].
\EE

For odd natural numbers $n$ and integers $m$ with $(m,n)=1$, the {\em Jacobi symbol}
\BD
   \JS mn
\ED
generalizes the Legendre symbol in the usual way (see \cite{Hu}, p. 44).
Note that the Jacobi symbol equals $1$ in the case $m=0$, $n=1$. If $n$ is even and
$(m, n)=1$, we put
\BD
   \JS mn=*,
\ED
where $*$ stands for an arbitrarily chosen symbol different from $\pm 1$. This means that the sequence of
convergents of $x$ defines a sequence
\BD
  \JS{s_k}{t_k},\: k\ge 0,
\ED
of Jacobi symbols. We call this sequence the {\em Jacobi sequence} of $x$ (although this name
is already in use in other fields of mathematics).
One of our main results is

\begin{theorem} 
\label{t1}

Let $e=2.718281\LD$ be Eulers's number. The Jacobi sequence of $e$ is purely-periodic
with period length 24. If, therefore $s_k/t_k$ is the $k$th convergent of $e$, then
\BD
   \JS{s_k}{t_k}=\JS{s_{k+24}}{t_{k+24}} \MB{ for all } k\ge 0.
\ED
Moreover, 24 is the smallest possible period length of the Jacobi sequence of $e$.

\end{theorem} 

\MN
It is easy to check that the period of the Jacobi sequence of $e$ reads
\BD
  1,1,-1,*,-1,*-1,-1,-1,*-1,*,\|-1,-1,1,*,1,*,1,1,1,*,1,*.
\ED
The symbol $\|$ separates the first half of the period from the second.
The latter arises from the former by interchanging $1$ and $-1$ (so one may say that
the period is skew-symmetric). Hence the period does not arise from a sub-period of length $12$,
but also not from one of length $8$.
Accordingly, 24 is the smallest possible period length.

\MN
A basic tool for our investigation (and, in particular, for the proof of Theorem \ref{t1}) is

\begin{theorem} 
\label{t2}

Let $x=[a_0,a_1,a_2,\LD]\in \R\smallsetminus\Q$ and $s_k/t_k$, $k\ge 0$, be as above. The Jacobi symbol $\JS{s_k}{t_k}$
depends only on the residue classes $\OV{a_0},\OV{a_1},\LD,\OV{a_k}\in\Z/4\Z$.
The same is true for the reciprocal symbol $\JS{t_k}{s_k}$.

\end{theorem} 

\MN
Possibly this theorem has been known implicitly, but we cannot give a reference where it is stated in the present
form.

In view of Theorem \ref{t2} we say that two irrational numbers $x=[a_0,a_1,a_2,\LD]$ and $y=[b_0,b_1,b_2,\LD]$
are {\em congruent mod 4}, if $\CMOD{a_k}{b_k}4$ for all $k\ge 0$. In this case we write $\CMOD xy4$.
Whenever $\CMOD xy4$, the numbers $x$ and $y$ have the same
Jacobi sequence. Of course, for any $x$ of the above shape there is a uniquely determined
$y=[b_0,b_1,b_2,\LD]$ such that
$\CMOD xy4$ and $b_k\in\{1,2,3,4\}$ for all $i\ge 0$.
We call $y$ the 4-{\em representative} of $x$ or, if we disregard $x$, a 4-representative {\em per se}.
The set of all possible 4-representatives has Lebesgue measure
0, since this is true for continued fractions with digits $a_i\le C$, $i\ge 1$, for an arbitrary constant $C$
(see \cite{RoSz}, p. 138). Hence congruence mod 4 divides $\R\smallsetminus\Q$ into a number of classes which can be
represented by a set of measure 0.

In the case of Euler's number $e$ we have
$e=[2,\{1,2j,1\}_{j=1}^{\infty}]$, where
$\{1,2j,1\}_{j=1}^{\infty}$ stands for the sequence
\BD
   1,2,1,1,4,1,1,6,1,\LD
\ED
(see \cite{Pe}, p. 124).
Therefore, the 4-representative of $e$ is $e'=[2,\{1,2,1,1,4,1\}_{j=1}^{\infty}]$, i. e., a periodic continued
fraction with period $1,2,1,1,4,1$, whose value is $e'=(7+\sqrt{15})/4$. In \cite{Gi2} we have shown that
the Jacobi sequence of a periodic continued fraction $x$ is periodic. More precisely, if the corresponding
{\em purely} periodic
continued fraction $z$ (here $z=[\{1,2,1,1,4,1\}_{j=1}^{\infty}]$) has a Jacobi sequence with even period length $L$,
then $L$ is also a possible period length for the Jacobi sequence of $x$. In \cite{Gi} we have shown that
a purely periodic continued fraction with period length $l$ has a periodic Jacobi sequence
with period length $L=dl$, where $d$ is a divisor
of $8$ or $12$. In our example $e'=[2,\{1,2,1,1,4,1\}_{j=1}^{\infty}]$ this means that $L$ can be chosen as a divisor
of $48$ or $72$. We shall show that $L=24$ works for $e'$ and, by Theorem \ref{t2}, also for $e$.

Theorem \ref{t2} says that for two irrationals $x$, $y$ with $\CMOD xy4$ the Jacobi sequences are the same.
Ist the converse also true, i. e., does equality of Jacobi sequences imply congruence mod 4? The answer is \BQ no\EQP,
as the following theorem shows.

\begin{theorem} 
\label{t3}

For each $C>0$ there is periodic 4-representative with period length $\ge C$ and the Jacobi sequence
$\{1\}_{j=1}^{\infty}=1,1,1,\LD$
Furthermore, there are uncountably many non-periodic 4-representatives having this Jacobi sequence.

\end{theorem} 

Below we shall investigate the Jacobi sequence for numbers like $e^2$, for which it is periodic with period length 40.
We shall also make the first steps towards a general theory of Jacobi sequences.
In particular, we study questions that arise in the context of Theorem \ref{t3}, such as: Does every possible period
length actually occur for some Jacobi sequence? Are there non-periodic Jacobi sequences? Are there short
sequences of symbols $\pm 1$ that do not occur as subsequences of Jacobi sequences --- for non-obvious reasons?

\section*{1. Proof of Theorem \ref{t2}}
Let the above notations hold, in particular $x=[a_0, a_1,a_2,\LD]$ (see (\ref{0.3})) and $s_k/t_k$ is defined by
(\ref{0.2}).
We have to show that the symbols $\JS{s_k}{t_k}$ and $\JS{t_k}{s_k}$ depend only on $\OV{a_0}, \LD, \OV{a_k}\in \Z/4\Z$.
The proof is by induction over $k$.
In the case $k=0$ we have $s_0=a_0, t_0=1$, and so $\JS{s_0}{t_0}=1$
and
\BD
 \JS{t_0}{s_0}=\left\{\begin{array}{cc}
                        1 & \MB{if } a_0 \MB{ is odd}, \\
                        * & \MB{otherwise.}
                      \end{array}\right.
\ED
For the step from $k$ to $k+1$ we put  $s=s_k$, $t=t_k$, $p=a_{k+1}$, $q=1$, $m=s_{k+1}$, $n=t_{k+1}$
and apply three theorems of \cite{Gi}. We have to distinguish a number of cases.

\MN
{\em Case 1}: $n=t_{k+1}$ is odd.

   (a) Let $t=t_k$ be odd. Then Theorem 1 of \cite{Gi} can be applied, since $q=1$ is odd.
It gives
\BD
  \JS{-\delta s}{t}\JS{p}{q}\JS{\delta m}{n}=\EPS(t,q,n).
\ED
Here the symbol $\EPS(t,q,n)$ equals 1, if at least two of $t,q,n$ are $\equiv 1$ mod 4, and $-1$, otherwise.
Since $q=1$, this symbol depends only on $\OV t$, $\OV n\in\Z/4\Z$, which, in turn, depend only on
$\OV{a_0}, \LD, \OV{a_{k+1}}\in \Z/4\Z$, by (\ref{0.2}). On the left hand side we have $\delta=(-1)^k$,
hence the Jacobi symbols $\JS{-\delta}{t}$, $\JS{\delta}{n}$ depend only on $k$ and $\OV t$, $\OV n\in\Z/4\Z$,
i. e., on $\OV{a_0}, \LD, \OV{a_{k+1}}\in \Z/4\Z$.
Altogether, we obtain, since $\JS pq=1$,
\BD
  \JS{m}{n}=\JS{-\delta}{t}\JS{\delta}{n}\EPS(t,q,n)\JS st.
\ED
By assumption, $\JS st$ depends only on $\OV{a_0}, \LD, \OV{a_{k}}\in \Z/4\Z$. As we have seen,
the remaining quantities on the right hand side depend only on $\OV{a_0}, \LD, \OV{a_{k+1}}\in \Z/4\Z$,
which, thus, holds for $\JS mn$.

As to the reciprocal symbol $\JS nm$, we have $\JS nm=*$, if $m$ is even,
so it depends only on $\OV m\in\Z/4\Z$ and thus on $\OV{a_0}, \LD, \OV{a_{k+1}}\in \Z/4\Z$.
If, however, $m$ is odd, we have, by quadratic reciprocity,
\BD
  \JS nm=\JS mn\EPS(m,n)
\ED
where $\EPS(m,n)=1$ if $m$ or $n$ is $\equiv 1$ mod 4, and $\EPS(m,n)=-1$, otherwise.
Again, $\JS nm$ depends only on $\OV{a_0}, \LD, \OV{a_{k+1}}\in \Z/4\Z$.

 (b) Let $t=t_k$ be even. Here we can apply Theorem 2 of \cite{Gi}. If $m$ is odd, this
theorem says
\BE
\label{1.1}
  \JS{\delta t}{s}\JS{p}{q}\JS{-\delta n}{m}=\EPS(s,q,m)
\EE
(observe that s is odd). Analogous considerations as in subcase (a) show that $\JS nm$ depends only
$\OV{a_0}, \LD, \OV{a_{k+1}}\in \Z/4\Z$. Further,
\BD
 \JS mn=\JS nm\EPS(n,m),
\ED
which gives the corresponding assertion for $\JS mn$.
If $m$ is even, the symbol $\JS nm (=*)$ depends only on $\OV m\in\Z/4\Z$ and, hence only, on
$\OV{a_0}, \LD, \OV{a_{k+1}}\in \Z/4\Z$.
As to $\rule{0mm}{5mm}\JS mn$, we have, by Theorem 2 of \cite{Gi},
\BE
\label{1.2}
  \JS{-\delta s}{s+t}\JS{p}{q}\JS{\delta m}{m+n}=\EPS(s+t,q,m+n).
\EE
Since $t$ is even, $s$ is odd, and, by quadratic reciprocity,
\BD
  \JS{s}{s+t}=\JS{s+t}{s}\EPS(s,s+t)=\JS ts\EPS(s,s+t),
\ED
so $\JS{s}{s+t}$ depends only on $\OV{a_0}, \LD, \OV{a_{k}}\in \Z/4\Z$, because this is true for $\JS ts$.
Therefore, the identity (\ref{1.2}) shows that $\JS{m}{m+n}$ depends only on $\OV{a_0}, \LD, \OV{a_{k+1}}\in \Z/4\Z$.
Since $n$ is odd and $\CMOD m{-n}{m+n}$, we have
\BD
\label{1.4}
    \JS{m}{m+n}=\JS{-n}{m+n}=\JS{-1}{m+n}\JS{m+n}{n}\EPS(n,m+n)=\JS{-1}{m+n}\JS{m}{n}\EPS(n,m+n),
\ED
where we have used quadratic reciprocity again. But $\JS{m}{m+n}$ depends only on $\OV{a_0}, \LD, \OV{a_{k+1}}\in \Z/4\Z$,
so the same must be true for $\JS mn$.

\MN
{\em Case 2}: $n=t_{k+1}$ is even.

Then $\JS mn = *$, so it depends only on $\OV n\in\Z/4\Z$ and, thus,
only on $\OV{a_0}, \LD, \OV{a_{k+1}}\in \Z/4\Z$. In the case of $\JS nm$ we observe that $t$ must be odd,
since $t_k, t_{k+1}$ cannot both be even. We apply Theorem 5 of \cite{Gi}.

(a) Suppose that $s$ is odd. Then this theorem says that (\ref{1.1}) holds.
By assumption, $\JS ts$ depends only on $\OV{a_0}, \LD, \OV{a_{k}}\in \Z/4\Z$. As in Part (b) of Case 1,
(\ref{1.1}) shows that $\JS nm$ depends only on $\OV{a_0}, \LD, \OV{a_{k+1}}\in \Z/4\Z$.

(b) Let $s$ be even. By the said theorem, (\ref{1.2})
holds in this case.
We use
\BD
  \JS{s}{s+t}=\JS{-t}{s+t}=\JS{-1}{s+t}\JS{s+t}{t}\EPS(s,s+t)=\JS{-1}{s+t}\JS st\EPS(s,s+t),
\ED
which shows that $\JS{s}{s+t}$ depends only on $\OV{a_0}, \LD, \OV{a_{k}}\in \Z/4\Z$.
Since $n$ is even, $m$ must be odd, and
\BD
  \JS{m}{m+n}=\JS{m+n}{m}\EPS(m,m+n)=\JS nm\EPS(m,m+n).
\ED
Together with (\ref{1.2}), this identity shows that $\JS nm$ depends only on
$\OV{a_0}, \LD, \OV{a_{k+1}}\in \Z/4\Z$.\\ \rule{5mm}{0mm} \STOP

\MN
{\em Remark.} It would be desirable to have a more elegant proof of Theorem \ref{t2},
in particular, a proof that avoids the above cases.

\section*{2. Jacobi sequences for $e$ and its relatives}

{\em Proof of Theorem 1.}
We start with the 4-representative $e'=[2,\{1,2,1,1,4,1\}_{j=1}^{\infty}]$ of Euler's number $e$.
The purely periodic number that belongs to $e'$ is $z=[\{1,2,1,1,4,1\}_{j=1}^{\infty}]$.
Let $s_k/t_k$ be the convergents of $z$ and $L$ an even multiple of the period length 6 of $z$.
Suppose that $L$ has the property
\BE
\label{1.6}
   \left(
     \begin{array}{cc}
       s_{L-1} & s_{L-2} \\
       t_{L-1} & t_{L-2} \\
     \end{array}
   \right)
   \equiv I \MB{ mod }4,
\EE
where $I$ is the $2\times 2$-unit matrix and the congruence has to be understood entry-by-entry.
Suppose, further, that
\BE
\label{1.8}
  \JS{t_{L-1}}{s_{L-1}}=1.
\EE
Then Proposition 1 of \cite{Gi} says
\BE
\label{1.10}
  \JS{s_k}{t_k}=\JS{s_{k+L}}{t_{k+L}}
\EE
for all $k\ge 0$.
If we choose $L=24$, we obtain (observe $z=(2+2\sqrt{15})/7)$
\BD
  \left(
     \begin{array}{cc}
       s_{23} & s_{22} \\
       t_{23} & t_{22} \\
     \end{array}
   \right)=
  \left(
     \begin{array}{cc}
       9286113 & 7622528\\
       6669712 & 5474849 \\
     \end{array}
   \right),
\ED
which is obviously $\equiv I$ mod 4. Moreover,
\BD
  \JS{t_{23}}{s_{23}}=\JS{6669712}{9286113}=1.
\ED
Hence (\ref{1.10}) holds for $L=24$.
The number $e'$ is mixed periodic, the number 2 forming its pre-period. We denote the sequence of
its convergents by
\BD
   \frac{p_0}{q_0}, \frac{s_0'}{t_0'}, \frac{s_1'}{t_1'},\frac{s_2'}{t_2'}, \LD,
\ED
in accordance with the pre-period of length 1. In \cite{Gi2} we have shown that our assumptions on $L$
imply
\BD
  \JS{s_k'}{t_k'}=\JS{s_{k+L}'}{t_{k+L}'}
\ED
for all $k\ge 0$.
In order to prove Theorem \ref{t1}, we have only to compare the Jacobi symbols $\JS{p_0}{q_0}$ and $\JS{s_{23}'}{t_{23}'}$.
Since they have the same value $(=1)$, the Jacobi sequence of $e'$ is purely periodic with period length 24.
\STOP

\MN
Next we consider the numbers $e^{1/n}$ for positive integers $n\ge 2$. We obtain

\begin{theorem} 
\label{t4}
The Jacobi sequence of $e^{1/n}$, $n\ge 2$, is purely periodic, the smallest possible period length being
\BD
   \left\{\begin{array}{rl}
            24 & \MB{ if } n \equiv 1,3 \MB{ mod } 4, \\
            12 & \MB{ if } n \equiv 2 \MB{ mod } 4,\\
            3 & \MB{ if } n \equiv 0 \MB{ mod } 4.
          \end{array}\right.
\ED
\end{theorem} 

\MN
{\em Proof.} By \cite{Pe}, p. 124, $e^{1/n}=[\{1, n(2j-1)-1,1\}_{j=1}^{\infty}]$, $n\ge 2$.
Accordingly, the 4-representative of $e^{1/n}$ is\\ \rule{15mm}{0mm}
\begin{tabular}{lll}
$[\{1,4,1,1,2,1\}_{j=1}^{\infty}]$ \rule{0mm}{5mm} &if&  $\CMOD n14$,\\
$[\{1,2,1,1,4,1\}_{j=1}^{\infty}]$ \rule{0mm}{5mm}&if& $\CMOD n34$,\\
$[\{1\}_{j=1}^{\infty}]$ \rule{0mm}{5mm}&if& $\CMOD n24$,\\
$[\{1,3,1,\}_{j=1}^{\infty}]$ \rule{0mm}{5mm} &if& $\CMOD n04$.
\end{tabular}

\MN
We investigate the period length of the Jacobi sequences of these four periodic continued fractions
in the same way as for $e'$, in particular, we use (\ref{1.6}) and (\ref{1.8}) for $L=24$ if $n\equiv 1,3$
mod 4, for $L=12$ if $\CMOD n24$, and for $L=6$ if $\CMOD n04$. In the last-mentioned case it turns out that
the smallest possible period length is not 6 but 3. Note that these cases are simpler than the case of $e$ since no
pre-period occurs.
\STOP

\MN
Even simpler than the case of $e^{1/n}$ is the case of the number $(e^{2/n}+1)/(e^{2/n}-1)$, $n\ge 1.$
Indeed, we have $(e^{2/n}+1)/(e^{2/n}-1)=[n,3n,5n,7n,\LD]\equiv[\{n,3n\}_{j=1}^{\infty}]$ mod 4 (see \cite{Pe}, p. 124).
Accordingly, the 4-representative of this number is
\\ \rule{15mm}{0mm}
\begin{tabular}{lll}
$[\{1,3\}_{j=1}^{\infty}]$ \rule{0mm}{5mm} &if&  $\CMOD n14$,\\
$[\{3,1\}_{j=1}^{\infty}]$ \rule{0mm}{5mm}&if& $\CMOD n34$,\\
$[\{2\}_{j=1}^{\infty}]$ \rule{0mm}{5mm}&if& $\CMOD n24$,\\
$[\{4\}_{j=1}^{\infty}]$ \rule{0mm}{5mm} &if& $\CMOD n04$.
\end{tabular}

\MN
If we inspect these cases in the above way, we obtain

\begin{theorem} 
\label{t5}
The Jacobi sequence of $(e^{1/n}+1)/(e^{1/n}-1)$, $n\ge 2$, is purely periodic, the smallest possible
period length being
\BD
   \left\{\begin{array}{rl}
            24 & \MB{ if } n \equiv 1,3 \MB{ mod } 4, \\
            8 & \MB{ if } n \equiv 2 \MB{ mod } 4,\\
            2 & \MB{ if } n \equiv 0 \MB{ mod } 4.
          \end{array}\right.
\ED
\end{theorem} 

\MN
Finally, we consider the number $e^2=[7,\{2+3(j-1),1,1,3+3(j-1),18+12(j-1)\}_{j=1}^{\infty}]$ (see \cite{Pe}, p. 125).
Its 4-representative is
\BD
e''=[3,\{2,1,1,3,2,1,1,1,2,2,4,1,1,1,2,3,1,1,4,2\}_{j=1}^{\infty}],
\ED
i. e.,
the quadratic irrational
\BD
  e''=\frac{370619+11\sqrt{444198255}}{177718}.
\ED
As in the case of $e$ we consider the purely periodic part of $e''$.
For this purely periodic number and $L=40$, (\ref{1.6}) holds, namely
\BD
  \left(
     \begin{array}{cc}
       11702972599281  & 5273785915232 \\
       4563573565840 & 2056512547601 \\
     \end{array}
   \right)
   \equiv I \MB{ mod }4.
\ED
Furthermore, (\ref{1.8}) is also fulfilled.
By the same arguments as in the case of the number $e$, we obtain
that $40$ is a possible period length of the Jacobi sequence of $e^2$.
Indeed, it turns out that this sequence is purely periodic with the period
\BEA
&& \rule{2mm}{0mm}1,*,\rule{2mm}{0mm}1,-1,*,\rule{2mm}{0mm}1,-1,*,-1,*,\rule{2mm}{0mm}1,* ,-1,-1,*,-1,-1,*,-1,*\\
&& -1,*,-1,\rule{2mm}{0mm}1,*,-1,\rule{2mm}{0mm}1,*,\rule{2mm}{0mm}1,*,-1,* ,\rule{2mm}{0mm}1,\rule{2mm}{0mm}1,*,\rule{2mm}{0mm}1,\rule{2mm}{0mm}1,*,\rule{2mm}{0mm}1,*.
\EEA
As in the case of $e$, we see that this period does not consist of subperiods of length 20 or 8.
Altogether, we obtain

\begin{theorem} 
\label{t6}

The Jacobi sequence of $e^2$ is purely-periodic
with smallest possible period length 40.

\end{theorem} 

\MN
{\em Remark.} It is not difficult to obtain analogues of Theorem \ref{t6} for the numbers
$e^{2/(2n+1)}$, $n\ge 1$ (their continued fraction expansion can be found in \cite{Pe}, p. 125). 
It turns out that the respective 4-representatives
have to be checked only for the cases $n=1,2,3,4$. In each of these cases one obtains 
40 as the smallest possible period length for
the Jacobi sequence. We leave these details to the reader.

\section*{3. Some period lengths of Jacobi sequences}

We start with a small table that displays the shortest possible period lengths of the Jacobi sequences of some periodic
continued fractions.

\bigskip\rule{2cm}{0cm}
\begin{tabular}{l|c} 

$[\{1\}_{j=1}^{\infty}]$, $[\{3\}_{j=1}^{\infty}]$ & 12\END
$[\{2\}_{j=1}^{\infty}]$& 8\END
$[\{4\}_{j=1}^{\infty}]$ & 2\END
$[1,1,\{4\}_{j=1}^{\infty}]$ & 1 \END
$[\{1,2,3\}_{j=1}^{\infty}]$ & 6\END
$[\{1,2,2\}_{j=1}^{\infty}]$ & 36 \END
$[\{1,2,2,2\}_{j=1}^{\infty}]$ & 8 \END
$[\{1,3,3\}_{j=1}^{\infty}]$ & 3

\end{tabular} 

\bigskip\noindent
This table shows that the shortest possible period length of a Jacobi sequence may be hard to predict from
the appearance of the continued fraction. But it leaves the impression that this period length is always a multiple of
the shortest possible period length of the latter. Theorem \ref{t3}, however, suggests that this
impression may be misleading. We are going to prove this theorem now.  A fundamental ingredient of the proof
is the following lemma.

\begin{lemma}  
\label{l1}
Let $x=[a_0,a_1,a_2,\LD]\in\R\smallsetminus\Q$ be such that the denominators $t_k$
and $t_{k+1}$, $k\ge 0$, are odd. 
Then
\BD
  \JS{s_k}{t_k}=\JS{s_{k+1}}{t_{k+1}},
\ED
except when
\BD
 \left\{ \begin{array}{l}
            \CMOD{t_k}14,\: \CMOD{t_{k+1}}34,\: k \MB{ odd}, \\
            \CMOD{t_k}34,\: \CMOD{t_{k+1}}14,\: k \MB{ even},\rule{0mm}{4mm}
          \end{array}\right.
\ED

\end{lemma} 
in which cases
\BD
  \JS{s_k}{t_k}=-\JS{s_{k+1}}{t_{k+1}}.
\ED

\MN
{\em Proof.} As in Section 1, we put $s=s_k$, $t=t_k$, $p=a_{k+1}$, $q=1$, $m=s_{k+1}$, $n=t_{k+1}$ and apply Theorem 1
of \cite{Gi}. In our case it says
\BD
  \JS{-\delta s}{t}\JS{\delta m}{n}=\EPS(t,n),
\ED
where $\delta=(-1)^k$. We have to distinguish eight cases depending on $t\equiv 1, 3$ mod 4,
$n\equiv 1, 3$ mod 4 and $k\equiv 0, 1$ mod 2. Thereby we obtain the lemma. \STOP

\MN
{\em Proof of Theorem \ref{t3}}.
Let $k_j$, $j\ge 1$, be a sequence of {\em even} natural numbers such that $k_1\ge 6$ and $k_{j+1}-k_j\ge 6$ for
all $j\ge 1$. We put
\BD
a_k=\left\{\begin{array}{ll}
             1 & \MB{if }k=0,1, \\
             2 & \MB{if } k=k_j-1 \MB{ or } k=k_j+1 \MB{ for some } j\ge 1, \\
             4 & \MB{ otherwise. }
           \end{array}\right.
\ED
Let $x=[a_0,a_1,a_2,\LD]$ and let $s_k/t_k$, $k\ge 0$, be the convergents of $x$.
Then
\BE
\label{3.2}
t_k\equiv\left\{\begin{array}{ll}
             3 \MB{ mod } 4 & \MB{if }k=k_j-1 \MB{ for some } j, \\
             1 \MB{ mod } 4  & \MB{ otherwise.}
           \end{array}\right.
\EE
In order to prove (\ref{3.2}), we put $k_0=0$ and use induction over $j$, $j\ge 0$.
Since $t_0=t_1=1$ and $t_k=4t_{k-1}+t_{k-2}$, we see that $\CMOD{t_k}14$ for $k_0+2\le k\le k_1-2$.
Suppose we have shown $\CMOD{t_k}14$ for $k_j+2\le k\le k_{j+1}-2$ and $j\ge 0$.
Then $t_{k_{j+1}-1}=2t_{k_{j+1}-2}+t_{k_{j+1}-3}$, and since $k_{j+1}-3\ge k_j+2$,
we have $t_{k_{j+1}-1}\equiv 2\cdot 1+1\equiv 3$ mod 4. Further,
$t_{k_{j+1}}=4t_{k_{j+1}-1}+t_{k_{j+1}-2}\equiv 4\cdot 3+1\equiv 1$ mod 4,
$t_{k_{j+1}+1}=2t_{k_{j+1}}+t_{k_{j+1}-1}\equiv 2\cdot 1+3\equiv 1$ mod 4,
and $t_{k_{j+1}+2}=4t_{k_{j+1}+1}+t_{k_{j+1}}\equiv 4\cdot 1+1\equiv 1$ mod 4.
Since $t_k=4t_{k-1}+t_{k-2}$, $k_{j+1}+3\le k\le k_{j+2}-2$, we have $\CMOD{t_k}14$ for these $k$.

Because of (\ref{3.2}) and Lemma \ref{l1},
$\JS{s_{k+1}}{t_{k+1}}=-\JS{s_k}{t_k}$ only if $k=k_j-2$ or $k=k_j-1$.
In the first case $\CMOD{t_k}14$ and $k$ is even, so the lemma says that this sign change is impossible.
In the second case $\CMOD{t_k}34$ and $k$ is odd, which excludes this sign change again. Hence the
symbol $\JS{s_k}{t_k}$ remains constant for all $k\ge 0$, and since $\JS{s_0}{t_0}=1$, the Jacobi sequence
of $x$ is $\{1\}_{j=1}^{\infty}$.
If we choose the numbers $k_j$ such that $k_{j+1}-k_j=k_1$ for all $j\ge 1$, the number $x$ is a periodic
4-representative with period length $k_1$, which can be made arbitrarily large. If we choose these numbers such
that the difference $k_{j+1}-k_j$ tends to infinity for $j\to \infty$, the 4-representative $x$ is not periodic.
Since there are uncountably many sequences $k_1,k_2,k_3,\LD$ with this property, we obtain uncountably many
4-representatives of this kind. \STOP

A modification of the construction of $a_k$ in the proof of Theorem \ref{t3} gives the following result.

\begin{theorem} 
\label{t7}

For each even number $L\ge 2$ there is a periodic 4-representative such that
its Jacobi sequence has the period $1,1,\LD,1,-1$, where the number $1$ is repeated $L-1$ times.
Furthermore, there are uncountably many non-periodic Jacobi sequences
belonging to 4-representatives $x\in\R\smallsetminus\Q$.

\end{theorem} 

\MN
{\em Proof.}
As in the proof of Theorem \ref{t3} let $k_j$, $j\ge 1$, be a sequence of {\em even} natural numbers
such that $k_1\ge 6$ and $k_{j+1}-k_j\ge 6$ for
all $j\ge 1$. We put
\BD
a_k=\left\{\begin{array}{ll}
             1 & \MB{if }k=0,1, \\
             2 & \MB{if } k=k_j \MB{ or } k=k_j+2 \MB{ for some } j\ge 1, \\
             4 & \MB{ otherwise. }
           \end{array}\right.
\ED
Let $x=[a_0,a_1,a_2,\LD]$ and let $s_k/t_k$, $k\ge 0$, be the convergents of $x$.
By the arguments of the proof of (\ref{3.2}),
\BD
\label{3.4}
t_k\equiv\left\{\begin{array}{ll}
             3 \MB{ mod } 4 & \MB{if }k=k_j \MB{ for some } j, \\
             1 \MB{ mod } 4  & \MB{ otherwise.}
           \end{array}\right.
\ED
Because of Lemma \ref{l1},
$\JS{s_{k+1}}{t_{k+1}}=-\JS{s_k}{t_k}$ only if $k=k_j-1$ or $k=k_j$.
In the former case we have $\CMOD{t_k}14$, $\CMOD{t_{k+1}}34$, and $k$ is odd, so this sign change actually occurs.
In the latter case we have  $\CMOD{t_k}34$, $\CMOD{t_{k+1}}14$, and $k$ is even, so the sign changes again.
Altogether, $\JS{s_k}{t_k}=1$ except for $k=k_j$, where $\JS{s_k}{t_k}=-1$.

As to the first assertion of the theorem, one chooses
the numbers $k_j$ such that $k_{j+1}-k_j=k_1$ for all $j\ge 1$. This proves the assertion for even numbers $L\ge 6$.
For $L=4$, we put $x=[1,1,4,\{4,2\}_{j=1}^{\infty}]$, which has the period $1,1,1,-1$,
and for $L=2$,
we put $x=[1,1,2,\{4\}_{j=1}^{\infty}]$, which has the period $1,-1$ (for all of these numbers $x$ the Jacobi sequence
has the pre-period 1).

As to the second assertion, we choose $k_j$ such that $k_{j+1}-k_j$ tends to infinity for $j\to \infty$. This gives
uncountably many non-periodic Jacobi sequences. \STOP

Theorem \ref{t7} shows that each even number $L\ge 2$ actually occurs as the period length of some Jacobi sequence.
Hence the case of odd numbers $L$ remains to be investigated.
The period $1,1,-1$ is impossible (see Theorem \ref{t9}). The construction
of the period $1,1,\LD,1,-1$ for the Jacobi sequence is always possible for odd period lengths $L\ge 5$,
but it is considerably more complicated than for even ones. Therefore, we restrict
ourselves to the following theorem which we obtain in a simpler way.

\begin{theorem} 
\label{t8}

For each natural number $L\ge 2$ there is a periodic 4-representative such that
its Jacobi sequence has the period $1,1,\LD,1,*$, where the number $1$ is repeated $L-1$ times.

\end{theorem} 

\MN
{\em Proof.} Let $L\ge 3$. We put
\BD
a_k=\left\{\begin{array}{ll}
             1 & \MB{if }\CMOD k{0,1}L, \\
             3 & \MB{if } \CMOD k{-1}L, \\
             4 & \MB{ otherwise. }
           \end{array}\right.
\ED
As above, we obtain
\BD
\label{3.64}
t_k\equiv\left\{\begin{array}{ll}
             0 \MB{ mod } 4 & \MB{if }\CMOD k{-1}L, \\
             1 \MB{ mod } 4  & \MB{ otherwise.}
           \end{array}\right.
\ED
Hence Lemma \ref{l1} shows that $\JS{s_k}{t_k}$ remains constant except in the cases $\CMOD k{-1}L$, where
$\JS{s_k}{t_k}=*$, and $\CMOD k0L$, where the lemma cannot be applied. Here, however,
we put $s/t=s_{k-2}/t_{k-2}$, $p/q=[3,1]=4/1$, $m/n=s_k/t_k$. Since $\JS pq=1$, Theorem 1 of \cite{Gi} yields
\BD
  \JS{-\delta s}t\JS{\delta m}n=\EPS(t,n).
\ED
Since $t\equiv n\equiv 1$ mod 4, we obtain $\JS mn=\JS st$, i. e., $\JS{s_k}{t_k}=\JS{s_{k-2}}{t_{k-2}}$.
Accordingly, $\JS{s_k}{t_k}=1$ if $k\not\equiv {-1}$ mod $L$, and $*$, otherwise.
In the case $L=2$, the Jacobi sequence of $[\{4\}_{j=1}^{\infty}]$ has the desired property. \STOP

\begin{theorem} 
\label{t9}

There is no number $x\in\R\smallsetminus\Q$ whose Jacobi sequence contains the subsequence
$-1,1,1,-1$. In particular, $1,1,-1$ cannot be the period of the Jacobi sequence of such
a number $x$.

\end{theorem} 

\MN{\em Proof.} Suppose $\JS{s_k}{t_k}=-1$, $\JS{s_{k+1}}{t_{k+1}}=\JS{s_{k+2}}{t_{k+2}}=1$,
$\JS{s_{k+3}}{t_{k+3}}=-1$.

{\em Case 1.} Let $k$ be odd. Then $k+2$ is also odd and due to the  sign change from $k+2$ to $k+3$,
$\CMOD {t_{k+2}}14$, 
by Lemma \ref{l1}. Since $k+1$ is even, $\CMOD {t_{k+1}}14$,
for otherwise there would be a sign change from $k+1$ to $k+2$, which is not the case.
Finally, $k$ is odd, and since there is a sign change from $k$ to $k+1$, we obtain $\CMOD {t_{k}}14$,
$\CMOD {t_{k+1}}34$. This, however, contradicts $\CMOD {t_{k+1}}14$.

{\em Case 2.} Let $k$ be even. Then $k+2$ is also even and the sign change from $k+2$ to $k+3$
requires $\CMOD {t_{k+2}}34$.
Now $k+1$ is odd, and since there is no sign change
from $k+1$ to $k+2$, we have $\CMOD {t_{k+1}}34$. But $k$ is even and there is a sign change from
$k$ to $k+1$, so $\CMOD {t_{k}}34$ and $\CMOD {t_{k+1}}14$. This is a contradiction again. \STOP

\MN
{\em Remark.} The proof also shows that $1,-1,-1,1$ is an impossible subsequence of a Jacobi sequence.


\vspace{0.5cm}
\noindent
Kurt Girstmair            \\
Institut f\"ur Mathematik \\
Universit\"at Innsbruck   \\
Technikerstr. 13/7        \\
A-6020 Innsbruck, Austria \\
Kurt.Girstmair@uibk.ac.at

\end{document}